\documentclass[a4paper,12pt]{article} 
\usepackage{amsmath,amsfonts,xcolor} 

\newtheorem{thm}{thm}[section]   
\newtheorem{theorem}[thm]{Theorem}   
\newtheorem{corollary}[thm]{Corollary}   
\newtheorem{proposition}[thm]{Proposition}   
\newtheorem{lemma}[thm]{Lemma}

\newcommand{\RR}{\mathbb{R}}   
\newcommand{\di}{\displaystyle}

\newcommand{\e}{\varepsilon}
\def\un{{\mathbf{1}}}

\begin{document}   
   
\title{\textbf{Sharp large time behaviour in $N$-dimensional reaction-diffusion equations of bistable type}}
\author{
{\bf Jean-Michel Roquejoffre}\\
Institut de Math\'ematiques de Toulouse; UMR 5219\\
Universit\'e de Toulouse; CNRS\\ 
Universit\'e Toulouse III,
118 route de Narbonne, 31062 Toulouse, France \\ 
\texttt{jean-michel.roquejoffre@math.univ-toulouse.fr}\\ 
\\[.5mm]
{\bf Violaine Roussier-Michon} \\ 
Institut de Math\'ematiques de Toulouse; UMR 5219\\
Universit\'e de Toulouse; CNRS\\ 
INSA Toulouse,
135 av. Rangueil, 31077 Toulouse, France\\ 
\texttt{roussier@insa-toulouse.fr}
}
\maketitle 
\begin{abstract}
We study the large time behaviour of the reaction-diffsuion equation $\partial_t u=\Delta u +f(u)$ in spatial dimension $N$, when the nonlinear term is bistable and the initial datum is compactly supported. We prove the existence of a Lipschitz function $s^\infty$ of the unit sphere, such that $u(t,x)$ converges uniformly in $\mathbb{R}^N$, as $t$ goes to infinity, to
 $U_{c_*}\bigg(|x|-c_*t +\di\frac{N-1}{c_*} \mathrm{ln}t + s^\infty\Big(\di\frac{x}{|x|}\Big)\bigg)$, where $U_{c*}$ is the unique 1D travelling profile. This extends  earlier results that identified the locations of the level sets of the solutions with $o_{t\to+\infty}(t)$ precision, or identified precisely the level sets locations for almost radial initial data.
\end{abstract}

\section{Introduction}    
\label{s1}

\subsection{Question under study}
The paper is devoted to the large time behaviour of the solution of the reaction-diffusion equation 
\begin{align}
\label{e1.1} \partial_t u= \Delta u +f(u), & \quad t>0 \, , \,  x \in \RR^N\\
\label{CI} u(0,x)=u_0(x), & \quad \quad \quad  \quad x \in \RR^N
\end{align}
where $f \in {\cal C}^\infty([0,1],\RR)$. We will assume the existence of $\theta\in(0,1)$ such that
$$
f<0\ \hbox{on $(0,\theta)$},\ f>0\ \hbox{on $(\theta,1)$},\quad f'(0)<0,\ f'(1)<0,\quad \int_0^1f>0.
$$
Thus  $f$ is said, in reference to the equation $\dot u=f(u)$, of the bistable type. A typical example is 
$$
f(u)=u(u-\theta)(1-u),\ \ \ 0<\theta<\frac12.
$$
We consider compactly supported  initial datum $u_0$  of the form
\begin{align}
\label{assumption}
\exists R_2>R_1>0 \, , \, \forall x \in \RR^N \, , \quad \un_{B_{R_1}}(x)\leq u_0(x) \leq \un_{B_{R_2}}(x),
\end{align}
where $\un_A$ is the indicator of the set $A$ and $B_{R}$ is the ball of $\RR^N$ 
of radius $R$ centered at the origin.
Equation \eqref{e1.1} has a unique classical solution $u(t,x)$ in $ {\cal C}^\infty([0,+\infty[ \times \RR^N, [0,1])$ emanating from $u_0$, see \cite{Hy} for instance. From Aronson and Weinberger \cite{aw}, as soon as $R_1>0$ is large enough, the solution $u$ spreads at a fixed speed $c_*>0$, in the following sense:
$$
\min_{|x|\leq ct} u(t,x) \to 1 \mbox{ as } t \to +\infty \, , \mbox{ for all } 0 \leq c < c_*
$$
and 
$$
\sup_{|x|\geq ct} u(t,x) \to 0 \mbox{ as } t \to +\infty \, , \mbox{ for all }  c > c_*.
$$
The goal of this work is to sharpen this result. The goal of this paper is to prove the following
\begin{theorem}
\label{thm1}
Let $u_0$ satisfy assumption \eqref{assumption}. There is a Lipschitz function $s^\infty$, defined on the unit sphere of $\RR^N$, such that the solution $u$ of \eqref{e1.1} emanating from $u_0$ satisfies
$$
\lim\limits_{t \to +\infty} \sup_{x \in \RR^N} \left|u(t,x) -U_{*}\biggl(\vert x\vert-c_*t+\frac{N-1}{c_*}{\mathrm{ln}}t+ s^\infty\Big(\frac{x}{\vert x\vert}\Big)\biggl)\right|=0.
$$
\end{theorem}

\subsection{Relation to existing works}

In the case $N=1$, equation \eqref{e1.1} with $N=1$ reads
\begin{equation}
\label{e1.11D}
\partial_t u=\partial_{xx} u +f(u),\quad t>0 \, , \, x\in\RR.
\end{equation}
It admits one-dimensional travelling fronts $U(x-ct)$ if and only if $c=c_*$, the just mentionned spreading speed. The profile~$U$, satisfies
\begin{equation}
\label{onde}
U'' +c_* \,U'+f(U)=0 , \quad x \in  \RR,
\end{equation}
together with the conditions at infinity 
\begin{equation} 
\label{CL onde}
\lim\limits_{x \to - \infty} U(x)=1 \quad\mbox{ and }\quad  \lim\limits_{x \to + \infty} U(x)=0. 
\end{equation} 
Any solution $U$ to \eqref{onde}-\eqref{CL onde} is a shift of a fixed profile $U_*$: $U(x)=U_*(x+s)$ with some fixed $s \in \RR$. The large time behaviour of  \eqref{e1.11D} has a history of important contributions, the most fundamental being perhaps that of Fife and McLeod \cite{FML}. They proved 
that the solution of~\eqref{e1.11D} starting from an initial datum $u_0(x)$ that is roughly front-like at infinity, namely
$$
\limsup_{x\to-\infty}u_0(x) > \theta,\quad\liminf_{x\to+\infty}u_0(x) < \theta
$$
gives rise to a solution $u(t,x)$ that converges to a travelling wave exponentially in time. Precisely, there exists $x_0\in\RR$ (depending on $u_0$ in a way that is not explicit in general) such that 
$$
\sup_{x\in\RR}\vert u(t,x)-U_*(x+c_*t+x_0)\vert\leq Ce^{-\omega t},
$$
where $\omega>0$ is essentially the first nonzero eigenvalue of the linear operator 
$$-\partial_{xx}-f'(U_*(.+x_0)).
$$

The large time behaviour of the solutions to \eqref{e1.1} has not been described at that level of precision, to the exception of a former paper by the second author \cite{VRM} tackling the case of almost spherically symmetric initial data, and which will be the starting point of our work.  This contribution proves the convergence to travelling waves, shifted by  the logarithmic delay $\di\frac{N-1}{c_*}\mathrm{ln} t$ plus an additional, possibly angle dependent constant. While the result takes advantage of the nearly spherically symmetric, it emphasises the fact that the part of the shift that is constant in time is, in general, angle dependent. We also refer to \cite{yagisita}, a work that also identifies the fact that alsmost spherically symmetric, but nonsymmetric initial data will remain so for all later times. 

In several space dimensions $N \geq 2$,  a line of results, in a spirit different from that of Theorem \ref{thm1}, is the convergence in profile of the solutions. Namely, $u(t,x)$  is followed in the reference frame where it is bounded away from 0 or 1, and its asymptotic shape is characterised. We mention a very interesting contribution of Jones \cite{J}, stating that the level sets of the solution of \eqref{e1.1}, whatever the nonlinearity is, will have oscillations only of the size $O_{t\to+\infty}(1)$. This is a consequence of the following fact: if $\lambda$ is a regular value of $u$, the normal to the $\lambda$-level set of $u$ meets the convex hull of the support of the initial datum. A  simple proof of this fact is given by Berestycki in \cite{B}. This work has been revisited in \cite{Ross1}.

Instead of a bistable nonlinearity, we may consider \eqref{e1.1} with a nonlinearity $f>0$ and concave between 0 and 1 (so-called Fisher-KPP nonlinearity, in reference to the seminal paper \cite{KPP}). It is a well-known fact that one-dimensional waves exist for all speeds $c\geq c_*=2\sqrt{f'(0)}$. If $U_*$ is a wave with bottom speed, we recently proved in \cite{RRR}, in collaboration with L. Rossi, that the dynamics of $u$ is
\begin{equation}
\label{e1.1200}
\lim\limits_{t \to +\infty} \sup_{x \in \RR^N} \left|u(t,x) -U_{*}\biggl(\vert x\vert-c_*t+\frac{N+2}{c_*}{\mathrm{ln}}t+ s^\infty\Big(\frac{x}{\vert x\vert}\Big)\biggl)\right|=0.
\end{equation}
Thus, in both cases, there is a logarithmic delay. However they are of different nature.   In the bistable case, the delay is purely due to curvature terms, as will be clear from Section 3, and as had already been elucidated in \cite{VRM}. In the Fisher-KPP case, there is an additional shift $\di\frac3{c_*}\mathrm{ln}t$, which is already present in one space dimension, and that is called the Bramson shift \cite{Br1}, \cite{Br2}. It comes from the fact that, as 0 is the most unstable value in the range of $u$ - that is, the growth for the linearised equation $\dot v=f'(u)v$ is maximal when $u=0$ -, the dynamics of $u(t,x)$ is driven by its tail, which implies a different behaviour that is very much related to the one-dimensional Dirichlet heat equation. Bramson's proof is probabilistic, and a new interpretation of this result is proposed in \cite{NRR}. Before the complete proof of \cite{RRR}, the position of the level sets had been identified with $O_{t\to+\infty}(1)$ precision by G\"artner \cite{G}, that is, they expand like $c_*t-\di\frac{N+2}2\mathrm{ln}t$.  Estimate \eqref{e1.1200} is proved on the basis of the ideas of \cite{RRR}.

As an illustration to our explanation, we mention the  recent contribution \cite{DQZ}, which  treats the porous medium equation with Fisher-KPP nonlinearity. It identifies the position of the level sets with $O_{t\to+\infty}(1)$ precision, that is, they expand like
$c_*t-\di\frac{N-1}{c_*}\mathrm{ln}t$. It may look surprising, as the nonlinearity is the Fisher-KPP one.  However, this can be explained by the fact that the porous medium equation is really a free boundary problem, so that the solution has no tail. This entails a  behaviour that is more closely related to what is observed in the bistable case.

\subsection{Strategy of the proof of Theorem \ref{thm1}, organisation of the paper}
Let us explain how the proof of Theorem \ref{thm1} proceeds. The first step is to identify the  reference frame in which $u(t,x)$ is nontrivial, for this we apply the existing analysis of the second author \cite{VRM}. Once this is done, we write, as in \cite{RRR}, equation \eqref{e1.1} in polar coordinates, shifted in the correct reference frame. This has the  inconvenience of cancelling out, at large times, the angular diffusion, which deprives us of an important source of compactness. To recover it we estimate the angular derivative, something that was quite useful in the fisher-KPP case \cite{RRR}. However, while one could use the maximum principle in \cite{RRR} in a relatively easy fashion - the asymptotic equation was the linear heat equation in the tail of the solution - one cannot do it here. Indeed, what drives the propagation is the body of the solution, not its tail. As a result, there is no obvious application of the maximum principle, and the estimate proceeds by applying a Fife-McLeod type idea to the angular derivative of $u$, by comparing it to its radial derivative. This is done in three successive steps detailed in Section 3. 
Once this is under control, a stability result, once again in the Fife-Mc Leod type, but complicated by the presence of angular terms, helps concluding the proof.

\noindent The organisation of the paper follows the main steps of this strategy.  In Section 2, we trap the solution between two 1D travelling waves moving like $c_*t-\di\frac{N-1}{c_*}\mathrm{ln}t$, thus characterising the reference frame in which the solution is nontrivial; we also prepare the equations. Section 3 is devoted to the main estimate, namely an estimate on the angular variable of $u$. The proof of Theorem \ref{thm1} is concluded in Section 4. We make some final remarks in Section 5.

\section{Radial bounds and preparation of the equations}
\label {section radial}

The main result of this section, that we will deduce from Theorem 1 of \cite{VRM}, is the following.
\begin{proposition}
\label{p2.1}
Let $u$ solve \eqref{e1.1} with initial datum $u_0$ satisfying \eqref{assumption}. There are four real numbers $t_0>0$, $C>0$ and $s_-<s_+$ such that, for all $t\geq t_0$ and $x \in \RR^N$,
we get
\begin{equation}
\label{e2.2}
U_*(|x|-c_*t + \frac{N-1}{c_*} \mathrm{ln}t -s_-)- C \, \frac{{\mathrm{ln}}t}{t}\leq u(t,x)\leq U_*(|x|-c_*t + \frac{N-1}{c_*} \mathrm{ln}t -s_+)+C \, \frac{{\mathrm{ln}}t}{t}.
\end{equation}
\end{proposition}
\noindent{\bf Proof.}
Let $u_0$ satisfy assumption \eqref{assumption} and $u$ be the unique solution to \eqref{e1.1} emanating from $u_0$. Define $R_0>0$ and $\delta_0>0$, depending on the non-linearity $f$, as in theorem 1 of \cite{VRM}. 

We first build a super-solution named $\bar{u}$ as follows. Choose $\varepsilon \in (0, \di\frac{\delta_0}{\sqrt{R_2+1}})$ and $\bar{R}>R_0$ such that
$$\forall x \in \RR^N \, , \quad u_0(x) \leq  \un_{B_{R_2}}(x) \leq U_*(|x| -\bar{R}) + \varepsilon \un_{B_{R_2+1}}(x)$$
Let $\bar{u}$ be the solution to \eqref{e1.1} emanating from $U_*(|x| -\bar{R}) + \varepsilon \un_{B_{R_2+1}}(x)$. By the maximum principle, we get that for all $t>0$ and all $x \in \RR^N$, $u(t,x) \leq \bar{u}(t,x)$ and we just have to compare $\bar{u}$ with a front. This is done by theorem 1 in \cite {VRM}. Indeed, defining 
$$X=\{u:\RR^N \to \RR \, |\, \exists \tilde{u}\in H^1(\RR^+) \mbox{ such that } u(x)=\tilde{u}(|x|) \mbox{ for } x \in \RR^N\}$$
we get
$$\|\bar{u}(0,x)-U_*(|x|-\bar{R})\|_X \leq \e \sqrt{R_2+1}\leq \delta_0$$
and therefore, by theorem 1 and the remarks below in \cite{VRM}, there exist $L \in \RR$ and $C>0$ such that for all $t>0$ and $x \in \RR^N$,
$$ | \bar{u}(t,x) - U_*(|x|-c_*t + \frac{N-1}{c_*} \ln t + L) |\leq C \, \frac{\ln t}{t}$$
Defining $s^+=-L$ leads to the right hand side of \eqref{e2.2} for any $t \geq 0$.

Dealing with a sub-solution is not so simple because a small perturbation as $max( U_*(|x|+\underline{R})+\e, 0)$ may not developp a front. We therefore use Aronson and Weinberger's result \cite{aw} to wait until the solution $u$ has propagated enough. Fix $\e>0$ and $\underline{R}>R_0$ such that $U_*(-\underline{R})\leq 1-\e$. Then, for all $x \in \RR^N$,  we have
$$U_*(|x|-\underline{R}) \leq U_*(-\underline{R}) \leq 1-\e.
$$

On the other hand, by Aronson and Weinberger's result \cite{aw} with $c=c_*/2$, there exists $t_ {\e} >0$ such that for any $t \geq t_\e$ and $|x|\leq ct$, $ 1-\e \leq u(t,x) \leq 1$. Choose $t_0 \geq t_\e$ such that 
$$\|U_*(\cdot -\underline{R})\|_{H^1(ct_0,\infty)} \leq \delta_0
$$ 
and define $\underline{u}(t_0,x)=U_*(|x|-\underline{R})\un_{B_{ct_0}}(x)$. Then, $\underline{u}(t_0,x) \leq u(t_0,x)$ for any $x \in \RR^N$.

Let $\underline{u}$ be the solution to \eqref{e1.1} emanating from $\underline{u}(t_0,x)$ at $t=t_0$. The maximum principle ensures that $\underline{u}(t,x) \leq u(t,x)$ for any $t\geq t_0$ and $x \in \RR^N$ and we just have to compare $\underline{u}$ with a front. Since $\|\underline{u}(t_0,x)-U_*(|x|-\underline{R})\|_X =\|U_*(\cdot -\underline{R})\|_{H^1(ct_0,\infty)} \leq \delta_0 $, theorem 1 in \cite {VRM} applies and there exists $L \in \RR$ and $C>0$ such that forall $t\geq t_0$ and $x \in \RR^N$, 
$$ | \underline{u}(t,x) - U_*(|x|-c_*t + \frac{N-1}{c_*} \ln t + L) |\leq C \, \frac{\ln t}{t}$$
Defining $s_-=-L$ leads to the left hand side of \eqref{e2.2} for any $t \geq t_0$. \hfill$\Box$

This proposition makes it clear that the transition zone, where $u$ is neither close to $1$ nor $0$, is located around $R(t)=c_*t-\di\frac{N-1}{c_*} \ln t$. We therefore choose to handle the initial equation \eqref{e1.1} in a frame, moving at speed $\dot{R}(t)$ in any radial direction. Let us explain those transformations on the equations.

From now on, we take $t=1$ as initial time and \eqref{CI} is replaced by $u(1,x)=u_0(x)$. This will be handier in view of the following transformations and, since equation \eqref{e1.1} is invariant by translation in time, there is no loss of generality.  We first use the polar coordinates 
$$
x\mapsto (r=\vert x\vert>0,\Theta=\frac{x}{\vert x\vert}\in \mathbb{S}^{N-1})
$$
then \eqref{e1.1} becomes
$$
\partial_t u=\partial_{rr}u +\frac{N-1}r\partial_ru +\frac{\Delta_\Theta u}{r^2}+ f(u),\quad{t>1,
\ r>0,\ \Theta\in \mathbb{S}^{N-1}}.
$$
Here, $\Delta_\Theta$ is the Laplace-Beltrami operator on the unit sphere of $\RR^N$. Its precise expression will not be needed in the sequel. The initial condition reads $u(1,r,\Theta)=u_0(r,\Theta)$.

Since we mentionned that the transition zone  is located around $R(t)=c_*t- k \ln t $ with $k=(N-1)/c_*$, we choose the change of variables $r'=r-R(t)$ and $u(t,r,\Theta)=u_1(t,r-R(t),\Theta)$. We drop the primes and indexes, and \eqref{e1.1} becomes
\begin{equation}
\label{e2.1} 
\partial_t u= \partial_{rr}u+c_*\partial_ru+\biggl(\frac{N-1}{r+c_*t-k{\mathrm{ln}}t}-\frac{k}t\biggl)\partial_r u+\frac{\Delta_\Theta u}{(r+c_*t-k{\mathrm{ln}}t)^2}+f(u).
\end{equation}
The equation is valid for $t>1$, $r>-2t+k{\mathrm{ln}}t$, and $\Theta\in \mathbb{S}^{N-1}$ and the initial condition becomes $u(1,r,\Theta)=u_0(r+c_*,\Theta)$.

To unravel the mechanisms at work, our first guess is that the term in $\Delta_\Theta v$ will not matter too much, because it decays like $t^{-2}$ (an integrable power of $t$), except in the zone $r\sim -c_*t$, where we know (for instance \cite{aw}) that $u(t,r,\Theta)$ goes to 1 as $t\to+\infty$. This confirms the information given by proposition \ref{p2.1} that the dynamics is like that of the one-dimensional equation. On the other hand, in the advection term, we have
$$
\di\frac{N-1}{r+c_*t-k{\mathrm{ln}}t}\sim_{t\to+\infty}\frac{N-1}{c_*t},
$$
except for extremely large $r$. This is nonintegrable in $t$,  but we  balance it with the $\di\frac{k}t$ term since we chose
\begin{equation}
\label{eq sur k}
k=\frac{N-1}{c_*}.
\end{equation}
This heuristics confirms that $R(t)=c_*t -\di\frac{N-1}{c_*} \ln t$ is the right moving frame to observe the large time dynamics of \eqref{e1.1}.
In the sequel, we will keep the notation $k$, keeping in mind that $k$ is given by formula \eqref{eq sur k}.  Also, from now on, we will only consider solutions of \eqref{e2.1}.

\section{Boundedness of the angular derivative}

This section is devoted to the following estimate
\begin{theorem}
\label{t3.1}
Let $u$ solve \eqref{e2.1} with initial datum $u_0(\cdot+c_*,\cdot)$. Then, there is $C>0$ such that 
\begin{equation}
\label{e3.1}
\forall t \geq 1 \, , \quad
\Vert\nabla_\Theta u(t,.,.)\Vert_{L^\infty((-c_*t/2,+\infty)\times\mathbb{S}^{N-1})}\leq C.
\end{equation}
\end{theorem}

The proof of this theorem \ref{t3.1} is based on a bootstrap argument. We will first prove that the quantity on the left handside is an $o(t)$, which will allow us to prove that it is an $O(t^\e)$ for all $\e>0$, which will in turn leads us to $O(1)$. The main idea is to adapt the construction by Fife and McLeod \cite{FML} of sub and super-solutions, but at the level of the linear equation. The main ingredient is that $\partial_ru$ becomes bounded away from 0 on every compact set, so that it may serve as a comparison function. And so, the main step (namely section \ref{section O(t)}) will consist in comparing $\vert\nabla_\Theta u\vert$ to a suitable multiple of $\partial_ru$, as it almost satisfies the same equation. The main body of the work will consist in quantifying what this innocent "almost" means. This idea of using the longitudinal  derivative of the solution as a comparison tool (as oppose to that of the wave, which has a long history dating back to Fife-McLeod) was first used in \cite{R1}, to prove the convergence to travelling waves in  cylindrical geometry.  

\medskip
\noindent{\bf Proof of theorem \ref{t3.1}.} So, let $u$ solve equation \eqref{e2.1} with datum $u_0(\cdot+c_*,\cdot)$.  

\subsection{The $o(t)$ estimate}

Let us perform the revert change of variables explained in the previous section to come back to $u$ solution to equation \eqref{e1.1}.
Pick any direction $\Theta$ on the unit sphere. We may, even if it means rotating, assume that $\Theta=0$, so that we are looking in the direction $Ox_1$. Let $x'=(x_2,...,x_{N})$ be the coordinates orthogonal to the direction $Ox_1$. Consider the sector
$$
\Sigma_t=\{x \in \RR^N \, |\, x_1>0,\ \frac{\vert x'\vert}{x_1}\leq\frac{1}{t^{3/4}}\},
$$
Notice that for $x \in \Sigma_t$, when $x_1\sim t$, we have $\vert x'\vert\leq t^{1/4}$. Write \eqref{e1.1} in $\Sigma_t$, in the reference frame moving like $c_*t-k{\mathrm{ln}}t$ in the direction $Ox_1$, it reads $X_1=x_1-c_*t +k \ln t$, $u(t,x)=u_1(t,X_1,x')=u_1(t,x_1-c_*t+k\ln t,x')$ so that dropping indexes,
$$
\partial_t u= \Delta u +\left(c_*-\frac{k}t \right) \partial_1 u +f(u).
$$
We also have, because $r=|(X_1+c_*t-k{\mathrm{ln}}t,x')|$ in $\Sigma_t-(c_*-k{\mathrm{ln}}t)e_1$:
$$
r-c_*t+k{\mathrm{ln}}t=X_1+o_{t\to+\infty}(1),\ \hbox{uniformly in $ \Sigma_t-(c_*-k{\mathrm{ln}}t)e_1$}
$$
Proposition \ref{p2.1} and parabolic regularity implies that the trajectories $(u(T+t,X_1,x'))_{T>0}$ are relatively compact in $C^2([-\tau,\tau]\times\RR\times[-M,M]^{N-1})$ for all $\tau>0$ and $M>0$. If $u_\infty(t,X_1,x')$ is a limiting trajectory we have for $(t,X_1,x')\in\RR^{N+1}$
\begin{align}
\label{e3.3}
&\partial_tu_\infty = \Delta u_\infty +\left(c_*-\di\frac{k}t \right)\partial_1 u_\infty +f(u_\infty)\\
&U_*(X_1-s_-)\leq u_\infty(t,X_1,x')\leq U_*(X_1-s_+).
\end{align}
From Theorem 1.1 of \cite{JMR-VRM} there is $s_\infty\in\RR$ such that  for $(t,X_1,x')\in\RR^{N+1}$
$$
u_\infty(t,X_1,x')=U_*(X_1-s_\infty).
$$
Parabolic regularity implies
$$
\lim_{t\to+\infty}\vert \nabla_{x'}u(t,X_1,x')\vert=0,\ \hbox{ uniformly in }(t,X_1)\in\RR_+\times\RR \mbox{ and }x' \mbox{ on every compact set.}
$$
Let us translate this result in the variables of equation \eqref{e2.1}. Because 
\begin{equation}
\nabla_\Theta u(t,r,0)=(r+c_*t-k{\mathrm{ln}}t)\nabla_{x'}u(t,X_1,0),
\end{equation}
we have the expected estimate on $\nabla_\Theta u$ for $\Theta=0$. Note that our argument is uniform in the direction considered, so that we have in the end, for $u$ solution to \eqref{e2.1}
$$
\lim_{t\to+\infty}\frac{\Vert\nabla_\Theta u(t,.,.)\Vert_{L^\infty((-c_*t/2,+\infty)\times\mathbb{S}^{N-1})}}{(r+c_*t-k{\mathrm{ln}}t)}=0.
$$
Parabolic regularity again implies the following corollary.
\begin{corollary}
\label{c3.1}
We have
$$
\lim_{t\to+\infty}\frac{\Vert\Delta_\Theta u(t,.,.)\Vert_{L^\infty((-c_*t/2,+\infty)\times\mathbb{S}^{N-1})}}{(r+c_*t-k{\mathrm{ln}}t)^2}=0.
$$
\end{corollary}
We also extract from the preceding argument an additional corollary. 
\begin{corollary}
\label{c3.2}
For every $M>0$, there is $T_M>0$ and $\delta_M>0$, the function $M\mapsto \delta_M$ decreasing, such that
$$
- \partial_ru(t,r,\Theta)\geq\delta_M\ \hbox{for $t\geq T_M$, $-M\leq r\leq M$, $\Theta\in{\mathbb{S}}^{N-1}$}.
$$
\end{corollary}

\subsection{The $O(t^\e)$ estimate} 
\label{section O(t)}

For $u$ solution to \eqref{e2.1}, denote
$$
V(t,r,\Theta)=-\partial_ru(t,r,\Theta). 
$$
The equation for $V$ is
\begin{equation}
\label{e3.4}
\biggl(\partial_t+L(t)-f'(u)\biggl) V=-\frac{N-1}{(r+c_*t-k{\mathrm{ln}}t)^2}V+\frac{2\Delta_\Theta u}{(r+c_*t-k{\mathrm{ln}}t)^3} ,
\end{equation}
the expression of $L(t)$ being
\begin{equation}
\label{e3.5}
L(t)=-\partial_{rr}-c_*\partial_r -\biggl(\frac{N-1}{r+c_*t-k{\mathrm{ln}}t}-\frac{k}t\biggl)\partial_r-\frac{\Delta_\Theta}{(r+c_*t-k{\mathrm{ln}}t)^2}.
\end{equation}
If $\Theta=(\theta_1,...,\theta_{N-1})$, we set
\begin{equation}
\label{e3.7}
u_i=\partial_{\theta_i}u, 
\end{equation}
we have
\begin{equation}
\label{e3.6}
\biggl( \partial_t+L(t)-f'(u) \biggl)u_i=0.
\end{equation}
A super-solution for \eqref{e3.6} is looked for under the form
\begin{equation}
\label{e3.8}
\overline v(t,r,\Theta)=\xi(t)V(t,r,\Theta)+q(t),
\end{equation}
that is, almost as in Fife-McLeod \cite{FML}, up to the fact that we work at the level of the linearised equation. 
In the sequel, we will always study $\overline v$ in the range $r\geq- c_*t/2$, the range $[-c_*t+k{\mathrm{ln}}t,-c_*t/2)$ being taken care of by Proposition \ref{p2.1}. Also, it is enough to construct $\overline v$ for $t$ large enough.

We shall now introduce new notations to explain how we will choose the functions $\xi(t)$ and $q(t)$. Denote $F=\|f'\|_\infty>0$ and 
$$\e(t)=\big\|\frac{N-1}{r+c_*t-k{\mathrm{ln}}t}V-\frac{2\Delta_\Theta u}{(r+c_*t-k{\mathrm{ln}}t)^2}\big\|_\infty \, .$$ 
Then, $\e(t)$ is a nonnegative function tending to $0$ at infinity, this last estimate comes from Corollary \ref{c3.1}. 
Pick $M>0$ and $\delta>0$ such that 
$$
f'(u(t,r,\Theta))\leq-\delta\ \hbox{if $\vert r\vert\geq M.$}
$$
By corollary \ref{c3.2}, there exists $T_M>0$ and $\delta_M>0$ such that $V \geq \delta_M$ for $t \geq T_M$, $|r|\leq M$ and $\Theta \in {\mathbb{S}}^{N-1}$. Moreover, if $V^-$ stands for the negative part of $V$, we also get from corollary \ref{c3.2} that 
 $$
 \lim_{t\to+\infty}\lim_{N\to\infty}\sup_{r\in(-c_*t/2,+\infty)\backslash(-N,N)}V^-(t,r,\Theta)=0.
 $$
Then, for any $\eta>0$, there exist $T>0$ and $N>0$ such that $V(t,r,\Theta) \geq 0$ for $r \in  [-N,N]$  and $|V^-(t,r, \Theta) |\leq \eta$ for $r \in (-c_* t/2,+\infty)\backslash(-N,N)$. Choose now $\eta >0$ small enough so that $0<\eta < \min(\delta_M, \di\frac{\delta \delta_M}{2F})$ and $N>M$. 

Now that all those preliminaries are given, let define $\xi$ and $q$ as the unique solutions to the ODE system for $t>T$
\begin{equation}
\label{e3.10}
\left\{\begin{array}{rll}
\dot q +\di\di\frac{\delta}4 q=&\di\frac{\e(t)\xi}{c_*t/2-k{\mathrm{ln}}t} \, , \\
\dot\xi=&\di\frac{\delta +F}{\delta_M+\eta}q \, ,
\end{array}\right.
\end{equation}
the initial data at $(q(T),\xi(T))$ being nonnegative and sufficiently large so that $\overline v(T,r,\Theta)\geq\vert u_i(T,r,\Theta)\vert$. Then, for any $t \geq T$, $q(t) \geq 0$ and $\xi(t) \geq 0$ and we shall prove that $\overline{v}$ defined by \eqref{e3.8} with $\xi$ and $q$ verifying \eqref{e3.10} is a supersolution to  \eqref{e3.6}. We have
$$
\begin{array}{rll}
&\biggl(\partial_t+L(t)-f'(u)\biggl)\overline v\\
=&\dot\xi(t)V+\xi(t)\biggl(\partial_tV+L(t)V-f'(u)V\biggl)+\dot q-f'(u)q\\
\geq&\dot\xi(t)V + \dot q-f'(u)q-\di\frac{\e(t)\xi}{r+c_*t-k{\mathrm{ln}}t},
\end{array}
$$
Consider first  the set where $\vert r\vert\geq M$. Since $\dot\xi\geq0$, $q\geq 0$, $|V^-|\leq \eta$ and $\eta < \min(\delta_M, \di\frac{\delta \delta_M}{2F})$, 
$$
\begin{array}{rll}
\biggl(\partial_t+L(t)-f'(u)\biggl)\overline v \geq & -\dot\xi \eta + \dot{q} +\delta q -\di\frac{\e(t)\xi}{c_*t/2-k{\mathrm{ln}}t},\\
=&- \eta \di\frac{\delta +F}{\delta_M+\eta} q+ \dot q + \delta q -\di\frac{\e(t)\xi}{c_*t/2-k{\mathrm{ln}}t}\\
\geq& \dot q + \di \frac{\delta}4 q-\di\frac{\e(t)\xi}{c_*t/2-k{\mathrm{ln}}t} =0,
\end{array}
$$
%
In the range $\vert r\vert\leq M$, the super-solution condition is true since
%
$$
\begin{array}{rll}
\biggl(\partial_t+L(t)-f'(u)\biggl)\overline v \geq & \dot\xi \delta_M + \dot{q} -F q -\di\frac{\e(t)\xi}{c_*t/2-k{\mathrm{ln}}t},\\
=& \delta_M \di\frac{\delta +F}{\delta_M+\eta} q+ \dot q -F q -\di\frac{\e(t)\xi}{c_*t/2-k{\mathrm{ln}}t}\\
\geq& \dot q + \di \frac{\delta}4 q-\di\frac{\e(t)\xi}{c_*t/2-k{\mathrm{ln}}t} =0,
\end{array}
$$
Finally, for any $r \geq -c_*t/2$, $\overline{v}$ is a supersolution to \eqref{e3.6}. 

Set $\e>0$. In order to prove $| \nabla_\Theta u | $ is a $O(t^\e)$, it suffices to study $(q(t), \xi(t))$ as time goes to infinity. 
The equation being linear, it is enough to study it with $(q(T),\xi(T))=(1,1)$. 
The first equation in \eqref{e3.10} gives
\begin{align*}
q(t)\leq & e^{-\delta(t-T)/4}+\int_T^t  e^{-\delta(t-s)/4} \di\frac{\e(s) \xi(s) }{c_* s/2-k \ln s}ds \, ,\\
 \leq &  e^{-\delta(t-T)/4}+ C \e \xi(t) \int_T^t  \di\frac{e^{-\delta(t-s)/4} }{1+s}ds \leq C \e \di \frac{\xi(t)}{1+t}
\end{align*}
with $C$ a universal constant. Indeed, $\dot \xi \geq 0$ and $T$ can be chosen large enough so that for any $s \geq T$, $0 \leq \e(s) \leq \e$.  We also have estimated $c_*t-k{\mathrm{ln}}t$ by a suitably small multiple of $1+t$. Plugging this result in the second equation of \eqref{e3.10}, we get for any $t \geq T$
$$
\dot\xi(t)\leq C \e \frac{\xi(t)}{1+t},
$$
with a universal constant $C$. Therefore, $\xi(t)\leq(1+t)^{C\e}$, which is precisely the desired estimate. The same will hold for $q$.

\subsection{Conclusion}

We start from
$$
\vert\nabla_\Theta u(t,r,\Theta)\vert\leq Ct^{1/20}.
$$
Set $u_{ij}=\partial_{\theta_i\theta_j}u$, we have
$$
\biggl(\partial_t+L(t)-f'(u) \biggl)u_{ij}=f''(u)u_iu_j=O(t^{1/10}).
$$
Parabolic regularity implies (this involves a rescaling of $\Theta$ by $t$, then scaling back):
$$
\vert u_{ij}(t,r,\Theta)\vert\leq Ct^{21/20}.
$$
This allows in turn a more precise estimate in the equation for $V=-\partial_r u$, because we have now
$$
\biggl(\partial_t+L(t)-f'(u)\biggl) V=O(t^{-39/20}),
$$
in the range $r\in(-c_*t/2,+\infty)$. Clearly, the right handside is an integrable power of $t$, we may therefore redo the previous step, replacing $\di\frac{\e(t)\xi(t)}{1+t}$ by$\di\frac{\xi(t)}{(1+t)^{39/20}}$. This leads to the desired estimate on $\overline v$ and then on $\nabla_\Theta u$, ending the proof of theorem \ref{t3.1}. \hfill$\Box$.

\section{Convergence}

Theorem \ref{thm1} will result from the following stability result, once again close in spirit to Fife and McLeod \cite{FML}.
\begin{theorem}
\label{t4.1}
Let $u$ be a solution to \eqref{e2.1} and $s$ a Lipschitz function of ${\mathbb{S}}^{N-1}$. For every $\e>0$, there is $T_\e>0$ and $\eta_\e>0$ (depending possibly on $\Vert s\Vert_\infty$ and $\Vert\nabla s\Vert_\infty$) such that, if we have
\begin{equation}
\label{hyp sur u}
\vert u(T_\e,r,\Theta)-U_*(r+s(\Theta))\vert\leq\eta_\e,\quad r>-c_*T_\e+{\mathrm{ln}}T_\e \, , \, \Theta \in {\mathbb{S}}^{N-1},
\end{equation}
then there is $T_\e'\geq T_\e$ such that, for all $t\geq T_\e'$ we have:
$$
\vert u(t,r,\Theta)-U_*(r+s(\Theta))\vert\leq\e,\quad r>-c_*T_\e'+{\mathrm{ln}}T_\e' \, , \, \Theta \in {\mathbb{S}}^{N-1}.
$$
\end{theorem}

Let us postpone the proof of theorem \ref{t4.1} to prove theorem \ref{thm1}.

\textbf{Proof of Theorem \ref{thm1}:}
Let $u$ solve \eqref{e1.1} with initial datum $u_0$ satisfying \eqref{assumption}. Perform transformations listed in section \ref{section radial} to deal with polar coordinates in the radial moving frame. We still denote $u$ the solution of \eqref{e2.1}. Define, for all $\tau>0$,
$$
\Omega_\tau=\{t\in[-\tau,\tau],r \in [-c_*t+k{\mathrm{ln}}t,+\infty),\Theta\in{\mathbb{S}}^{N-1}\}
$$
Then, by parabolic regularization and Ascoli's theorem, the family $(u(T+.,.,.))_{T>0}$ is relatively compact in $C(\Omega_\tau)$ for every $\tau>0$. Therefore, there is a sequence $(t_n)_n$ going to infinity such that $(u(t_n+.,.,.))_n$ converges, uniformly in every $\Omega_{\tau}$, to a uniformly continuous
function $u_\infty(t,r,\Theta)$ satisfying for $t \in \RR$
$$
\begin{array}{rll}
\partial_tu_\infty= \partial_{rr}u_\infty +c_* \partial_r u + f(u_\infty)\quad&\hbox{in ${\mathcal{D}}'(\RR^2\times{\mathbb{S}}^{N-1})$} \\
U_*(r-s_-)\leq u_\infty(t,r,\Theta)\leq U_*(r-s_+).&
\end{array}
$$
thanks to the radial barriers obtained in section \ref{section radial}.
Moreover, $u_\infty(.,.,\Theta)$ is  $C^1$ in $t$ and $C^2$ in $r$ due to parabolic regularity, for every $\Theta\in{\mathbb{S}}^{N-1}$. From \cite{FML}, for every $\Theta\in{\mathbb{S}}^{N-1}$, there is $s(\Theta)\in\RR$ such that
$$
u_\infty(t,r,\Theta)=U_*(r+s(\Theta)).
$$
From Theorem \ref{t3.1}, the function $s$ is Lipschitz. Moreover, from Theorem \ref{t4.1}, the whole family $(u(t,.,.))_{t>0}$ converges uniformly on $[0,+\infty) \times \mathbb{S}^{N-1}$ to the function $(r,\Theta)\mapsto U_*(r+s(\Theta))$. Reverting to the original variables proves theorem \ref{thm1}.  \hfill$\Box$

As said above, Theorem \ref{t4.1} is a Fife and Mc Leod's type result, and so will be obtained from the construction of sub and super solutions very much inspired from \cite{FML}. However, while those in \cite{FML} were explicit, the ones constructed here
solve a nonlinear differential system, and so must be studied with a little care. This is the object of the following intermediate lemma. 

For any $\e>0$, $C>0$ and $T_0>0$ such that $c_*T_0-k{\mathrm{ln}}T_0\geq100$, we define
$$
g_\e(t,\xi)=\frac{C}\e\big\vert\frac{N-1}{c_*t-k{\mathrm{ln}}t+\xi}-\frac{k}t\big\vert \, , \, t>T_0 \, , \,  \xi \in \RR
$$
\begin{lemma}
\label{l4.1}
Let $\delta$, $\gamma$ and $\eta$  be strictly positive constants. Consider the differential system
\begin{equation}
\label{e4.1}
\left\{
\begin{array}{rlll}
\dot q+\delta q=&g_\e(t,\xi) & t > T\\
\gamma\dot\xi=&Cq+g_\e(t,\xi) & t > T\\
q(T)=&\eta,\ \xi(T)=0. &
\end{array}
\right.
\end{equation}
where $\e$, $C$ and $g_\e$ are defined above. 
There is $T \geq T_0$ and $K>0$ depending on all the constants involved in \eqref{e4.1}, except $\eta$, such that \eqref{e4.1} has a solution defined on $[T,+\infty)$ which satisfies for any $t\geq T$, 
\begin{equation}
\label{e4.2}
0\leq q(t)\leq K(\eta+\frac1{\e T^{1/2}}) e^{-\delta/2(t-T)}+\frac{K}{t^{3/2}},\quad \xi(t)\leq K(\eta+\frac1{\e T^{1/2}}).
\end{equation}
\end{lemma}
\noindent{\bf Proof.}   In what follows, $K$ denotes a generic positive constant that may differ from lines to lines. We first derive a logarithmic bound on $\xi$, then the desired bound. By definition, there exists $K>0$ such that for any $t>T$, 
$$
0 \leq g_\e(t,\xi(t))\leq \frac{K}{\e(t+\xi(t))},
$$
so that we have
$$
q(t)\leq \eta e^{-\delta(t-T)}+\frac{K}\e\int_T^t\frac{e^{-\delta(t-s)}}{s+\xi(s)}ds.
$$
Let $T^* \geq T$ be the largest $\tau \geq T$ such that
\begin{equation}
\label{e4.3}
\hbox{for all $t\in[T,\tau]$,}\ \xi(t)\leq({\mathrm{ln}}t)^2.
\end{equation}
So, for any $s \in [T,T^*]$, we may write
$$
\frac1{\e(s+\xi(s))}\leq \frac{K}{\e s},
$$
so that for any $t \in [T,T^*]$, cutting the integral at $s=\di\frac{t}2$, we get
$$
q(t)\leq \frac{K}{\e t}+(\frac{K}{\e T}+\eta)e^{-\frac{\delta}{2}(t-T)}.
$$
This implies in turn
$$
\dot\xi(t)\leq \frac{K}{\e t},
$$
hence $\xi(t)\leq \di\frac{K}\e{\mathrm{ln}}t$, a contradiction with the definition of $T^*$, as soon as $T$ is large enough, say, of the order $\e^{-2}$. So for any $t \geq T$,  $\xi(t)\leq({\mathrm{ln}}t)^2$. But then, we have a more precise estimate of 
$g_\e(t,\xi(t))$. Using the actual value of $k$ we have
$$
g_\e(t,\xi(t))\leq \frac{C(N-1)}{\e c_*t}\biggl(\frac1{1-\ln t/t+({\mathrm{ln}}t)^2/(c_*t)}-1\biggl)\leq\frac{K{\mathrm{ln}}t}{\e t^2}\leq\frac{K}{\e t^{3/2}},
$$
thus an integrable power of $t$. Plugging this new estimate into \eqref{e4.1} yields \eqref{e4.2}, hence the lemma.  \hfill$\Box$

\noindent{\bf Proof of Theorem \ref{t4.1}.} 
Let $u$ be a solution to \eqref{e2.1} and $s$ as in the statment of theorem \ref{t4.1}. Let $\e>0$ and choose $T_\e =O(\frac{1}{\e^2})$ and $\eta_\e=O(\e^2)$. 
We first regularise $s$: set
$$
s^\pm_\e(\Theta)=(\rho_\e*s)(\Theta)\pm C\e,
$$
where $\rho_\e$ is an $\e$-approximation of the identity and $C>0$ large enough (say, a multiple of $\Vert\nabla s\Vert_\infty$). Then,  by assumption \eqref{hyp sur u}, we have
$$
U_*(r+s^-_\e(\Theta))-C\e\leq u(T_\e,r,\Theta)\leq U_*(r+s^+_\e(\Theta))+C\e.
$$
Also notice that 
$$
\Vert\nabla_\Theta s^+_\e\Vert_\infty\leq\Vert\nabla_\Theta s\Vert_\infty,\ \Vert\Delta_\Theta s^+_\e\Vert_\infty\leq\frac{\Vert\nabla_\Theta s\Vert_\infty}\e.
$$
Write equation \eqref{e2.1} for $u$ as 
$$
N\!L[u]=0.
$$
We are going to construct sub and super-solutions to \eqref{e2.1}, starting from $t=T_\e$, at small distance of $u(t,r,\Theta)$. We display details for the super-solution's case, as it goes in the same way for sub-solutions. We try the form, in the spirit of Fife and Mc Leod \cite{FML},
\begin{equation}
\label{e4.5}
\overline u(t,r,\Theta)=U_*(r+s^+_\e(\Theta)-\xi^+(t))+q^+(t).
\end{equation}
where $\xi^+$ and $q^+$ are chosen in the following way.

We first define $\mu_0\in(0,\min(\theta, 1-\theta)))$ and $\delta>0$ such that 
$$
f'(u)\leq-\delta\ \hbox{on }[0,\mu_0]\cup[1-\mu_0,1].
$$
Choose $M>0$ such that $U_*(\rho) \in [0,\mu_0]\cup[1-\mu_0,1]$ if $\rho \notin [-M,M]$ and $\delta_M>0$ such that $U'_*(\rho) \leq -\delta_M$ if $\rho \in [-M,M]$.

Now define $\xi^+$ and $q^+$ as the unique solution to \eqref{e4.1} with $\gamma =\delta_M$ and $\eta=C\e$. Then, choosing $T_\e$ possibly larger, we know by lemma \ref{l4.1} that there exists $K>0$ such that for any $t \geq T_\e$, 
$$
0\leq q(t)\leq K(C\e+\frac1{\e T^{3/2}}) e^{-\delta/2(t-T)}+\frac{K}{t^{3/2}},\quad \xi(t)\leq K(C \e+\frac1{\e T^{1/2}}).
$$
Adjusting $\e >0$ such that $2K(C\e+\frac{1}{\sqrt{T_\e}}) < \mu_0$, it implies that
$$0 \leq q^+(t) \leq \mu_0 \mbox{ and } \dot{\xi}^+(t) \geq 0$$

Once these preliminaries are done with, we may compute $NL[\bar{u}]$. Set $\rho:=r+s^+_\e(\Theta)-\xi^+(t)$. Then,
$$
\begin{array}{rll}
N\!L[\overline u]=&-\dot\xi^+(t) U_*'(\rho)+\dot{q}^+(t) -f(U_*(\rho)+q^+)+f(U_*(\rho))\\
&-\biggl(\di\frac{N-1}{\rho+c_*t-k{\mathrm{ln}}t+\xi^+(t)-s^+_\e(\Theta)}-\di\frac{k}t\biggl)U_*'(\rho)-\di\frac{U_*'(\rho)\Delta_\Theta s^+_\e+U_*''(\rho)\vert\nabla_\Theta s^+_\e\vert^2}{(\rho+c_*t-k{\mathrm{ln}}t+\xi^+(t)-s^+_\e(\Theta))^2}.
\end{array}
$$
Consider first the range $\vert\rho\vert\geq M$. Then, plugging $\overline u$ in the expression of $N\!L$ reveals that a sufficient condition for $N\!L[\overline u]\geq 0$ is
\begin{equation}
\label{e4.12}
\dot q^++\delta q^+\geq\frac{\Vert\nabla_\Theta s\Vert_\infty}{\e\biggl(\rho+c_*t-k{\mathrm{ln}}t+\xi^+(t)\biggl)^2}+\biggl\vert \frac{N-1}{\rho+c_*t-k{\mathrm{ln}}t+\xi^+(t)}-\frac{k}t\biggl\vert.
\end{equation}
In the range $\rho\leq M$, a sufficient condition is
\begin{equation}
\label{e4.14}
\delta_M\dot \xi^+-C q^+\geq\frac{\Vert\nabla_\Theta s\Vert_\infty}{\e\biggl(\rho+c_*t-k{\mathrm{ln}}t+\xi^+(t)\biggl)^2}+\biggl\vert \frac{N-1}{\rho+c_*t-k{\mathrm{ln}}t+\xi^+(t)}-\frac{k}t\biggl\vert,
\end{equation}
Both conditions are satisfied as $\xi^+$ and $q^+$ satisfy \eqref{e4.1}. Note that  we always may absorb the first term in the right handside of \eqref{e4.12} or \eqref{e4.14} in a function of the type $g_\e(t,\xi^+(t))$. Thus $\bar{u}$ is a supersolution to \eqref{e2.1} and for any $t \geq T_\e$
$$u(t,r,\theta) \leq U_*(r+s(\theta))+\e$$
Dealing in the same way with a subsolution finishes the proof of the theorem. \hfill$\Box$.

\section{Further questions and final remarks}

Instead of working with a bistable nonlinearity, one could think of dealing with a nonlinearity $f$ for which there is $\theta\in(0,1)$ such that
$$
f\equiv0\ \hbox{on $[0,\theta]$},\ \ \ f>0\ \hbox{on $(\theta,1)$},\ \ \ f(1)=0.
$$
We believe that Theorem \ref{thm1} holds in this case. However, the degeneracy of $f$ near 0 would certainly impose the construction of sub and super solutions with exponential weights.    The main step would once again be the gradient estimate, which already carries more technicalities than the basic Fife-Mc Leod construction. Some rather tedious computations should therefore be expected.

A more challenging question is to assess when the shift $s^\infty(\Theta)$ is trivial, that is, constant with respect to $\Theta$. The analysis of the second author in \cite{VRM} shows that the set of (almost spherically symmetric) initial data giving raise to a nontrivial shift is quite big (open and dense). The issue is therefore whether any non spherically symmetric initial datum will generate a nontrivial shift. We note that \cite{VRM} does not preclude codimension 1 sets of initial data generating trivial shifts. We hope to address this question in the future.

Another, possibly easier question, concerns further regularity for $s^\infty$. In \cite{VRM} it is shown to be $L^2$, in this work we upgrade its regularity to Lipschitz. Whether it is $C^2$, or whether the derivative may develop discontinuities even if the initial datum is smooth, is something that does not immediately stem from our analysis. We simply note that, given the goal that is ours in this paper, additional regularity is only anecdotical, in the sense that it would have slightly simplified the convergence proof. However the question is interesting in its own right and we note that, in \cite{RRR} we could not go further than Lipschitz regularity either. Whether more involved considerations would have allowed us to reach further regularity, or, as opposed to this, a new phenomenon occurs, is something that we do not know.
 
Finally, we believe that it would be  interesting to understand the rate of convergence of $u(t,x)$ to the shifted wave. In he Fisher-KPP case, it is a very interesting problem that was first raised by Ebert and Van Saarloos \cite{ES}, in one spatial dimension. They proposed a full expansion in terms of powers of $t^{-1/2}$ and found universal (that is, not depending on the initial datum and not even on the nonlinearity $f$) terms. The analysis of \cite{ES} was carried out in the formal style, and a first mathematically rigorous proof of the expansion up to the order $t^{-1+\delta}$ ($\delta>0$ is any small  positive number) was given in \cite{NRR2}. An analysis, still in the formal style, of \cite{BBD}, finds an expansion that is different from that of \cite{ES} at higher orders, in the sense that $t^{-1}\mathrm{ln}t$ terms pop up. The analysis of \cite{BBD} is confirmed in a mathematically rigorous way by Graham in \cite{Gra}. So, coming back to our model \eqref{e1.1}, we believe that pushing the expansion, perhaps to exponentially small terms, would be of interest. Making a more extensive use of the concept of radial waves, as defined in \cite{VRM}, could be a starting point.




{\footnotesize
 
}

\end{document}